\documentclass[11pt,a4paper,english,reqno]{amsart}
\usepackage{amsmath,amssymb,amsfonts,epsfig,mathrsfs}
\usepackage[T1]{fontenc}

\usepackage{color}
\usepackage{array}
\usepackage{amsthm}
\usepackage{amstext}
\usepackage{graphicx}
\usepackage{setspace}
\usepackage[margin=2.5cm]{geometry}
\usepackage{bbm}
\usepackage{color}
\usepackage{enumitem}
\allowdisplaybreaks[4]

\usepackage{pgfplots}

\usepackage{amscd,psfrag}
\usepackage{yhmath}
\usepackage[mathscr]{eucal}

\usepackage{slashed}

\makeatletter
\pdfpageheight\paperheight
\pdfpagewidth\paperwidth

\usepackage{mathrsfs}

\usepackage{epstopdf}
\usepackage{chngcntr}
\counterwithin{figure}{section}
\usepackage{mathrsfs}

\usepackage{indentfirst}	

\usepackage[normalem]{ulem}
\theoremstyle{plain}

\newtheorem{definition}{Definition}[section]
\newtheorem{theorem}[definition]{Theorem}
\newtheorem*{theorem*}{Theorem}

\newtheorem{remark}[definition]{Remark}

\newtheorem*{remark*}{Remark}
\newtheorem*{sideremark*}{Side Remark}\newtheorem*{mt*}{Main Theorem}

\newtheorem*{claim*}{Claim}
\newtheorem*{q*}{Question}

\newtheorem*{corollary*}{Corollary}
\newtheorem*{proposition*}{Proposition}

\newcommand{\R}{\mathbb{R}}

\newcommand{\na}{\nabla}

\newcommand{\p}{\partial}
\newcommand{\e}{\epsilon}
\newcommand{\emb}{\hookrightarrow}

\newcommand{\map}{\rightarrow}

\newcommand{\A}{{\mathscr{A}}}
\newcommand{\haus}{\mathcal{H}}

\newcommand{\ball}{{\mathbb{B}}}
\newcommand{\gr}{{\bf Gr}}

\allowdisplaybreaks[4]

\def\XXint#1#2#3{{\setbox0=\hbox{$#1{#2#3}{\int}$ }
\vcenter{\hbox{$#2#3$ }}\kern-.6\wd0}}

\numberwithin{equation}{section}
\numberwithin{figure}{section}

\title{A Note on Generic Transversality of Euclidean Submanifolds}

\author{Siran Li}
\address{Siran Li: Department of Mathematics, Rice University, MS 136
P.O. Box 1892, Houston, Texas, 77251-1892, USA; \, $\bullet$ \,  Department of Mathematics, McGill University, Burnside Hall, 805 Sherbrooke Street West, Montreal, Quebec, H3A 0B9, Canada.}

\email{\texttt{Siran.Li@rice.edu}}


\date{\today}

\begin{document}

\maketitle

\section{The Main Result}
In this short note, we establish a quantitative description of the genericity of transversality of $C^1$-submanifolds in $\R^n$:
\begin{theorem}\label{thm: main}
Let $\Sigma \subset \R^n$ be a $d$-dimensional $C^1$-embedded submanifold where $n \geq d+1$ and $d \geq 1$. Denote by
\begin{equation*}
\A(\Sigma) := \bigg\{ a \in \R^n: {\rm volume}\,\Big\{ p\in\Sigma : \p\ball(a, |a-p|) \text{ is not transverse to $\Sigma$ at $p$} \Big\} > 0  \bigg\}.
\end{equation*}
Then $\A(\Sigma)$ is contained in a countable union of $(n-d-1)$-dimensional affine planes.
\end{theorem}

The theorem states that $\A(\Sigma)$, the ``exceptional set'' for transversality, is small in the measure-theoretic sense. Indeed, it is $\haus^{n-d-1}$-rectifiable (see Federer \cite{f}). Here $\haus^k = $ Hausdorff measure of dimension $k$, and $\ball(a,\rho)$ denotes the Euclidean ball centred at $a$ of radius $\rho$ in $\R^n$. The volume is taken with respect to the metric induced by the embedding $\Sigma \emb \R^n$.

Let $f: \Sigma \map \R^n$ be a differentiable map and let $\hat{\Sigma}\subset\R^n$ be a $C^1$-submanifold. By definition (\cite{gp}, Guillemin--Pollack), $f$ is {\em transverse} to $\hat{\Sigma}$ if and only if for every $p \in f^{-1}(\hat{\Sigma})$ there holds
\begin{equation*}
{\rm image} \,(d_pf) + T_{f(p)}\hat{\Sigma} = T_{f(p)}\R^n.
\end{equation*}
Moreover, if $\Sigma$ and $\hat{\Sigma}$ are both embedded submanifolds of $\R^n$, then they are said to be {\em transverse} to each other if and only if the embedding $\iota: \Sigma \emb \R^n$ is transverse to $\hat{\Sigma}$. Via the identification of $\iota$ with the natural inclusion of $\Sigma$ as a subset of $\R^n$, the definition is equivalent to that
\begin{equation*}
T_p\Sigma + T_p\hat{\Sigma} = T_p \R^n\qquad \text{ for every } p \in \Sigma \cap \hat{\Sigma}.
\end{equation*}

By well-known results in differential topology {\it \`{a} la} R. Thom (\cite{thom}; also {\it cf.} \cite{gp}), the transversality of Riemannian submanifolds is a  generic property. Theorem \ref{thm: main} aims at providing one quantitative result in this direction.

The case $d=1$ of Theorem \ref{thm: main} is proved in $\S 6$, \cite{hr} by R. Hardt and T. Rivi\`{e}re. We generalise the arguments therein to prove for general $d$. In \cite{hr}, Hardt--Rivi\`{e}re used the result for $d=1$ to study the non-density of smooth maps in Sobolev spaces between Riemannian manifolds, and established deep  connections between harmonic maps and the minimal model programme.

Let us also mention the recent paper \cite{bc} by A. Bressan and G. Chen, as well as many subsequent developments, for applications of the generic transversality of curves to nonlinear partial differential equations (PDE). It is interesting to explore the possible applications of Theorem \ref{thm: main} in the study of generic properties for PDE solutions.

The remaining parts of the note consist of \S $\ref{sec: proof}$, in which we prove Theorem \ref{thm: main}, and \S $\ref{sharp}$, in which we demonstrate the sharpness of the theorem by constructing an explicit example.

\section{Proof of Theorem \ref{thm: main}}\label{sec: proof}

Let us first introduce several notations and reductions:
\begin{itemize}
\item
For $n$ fixed as in Theorem \ref{thm: main}, we write $\gr_i$ for the $i$-Grassmannian in $\R^n$, {\it i.e.}, the $i$-dimensional subspaces of $\R^n$, or the $i$-planes in $\R^n$ passing through the origin. 
\item
$pr_i:\R^n \times \gr_i \map \R^n$ is the projection onto the first variable, namely the projection of the $i$-Grassmannian bundle. 
\item
As $\Sigma \subset \R^n$ is a Riemannian manifold (second countable and para-compact), it suffices to prove Theorem \ref{thm: main} on each coordinate chart. 

So, from now on we may assume that $\Sigma = {\rm image}(\Phi)$ where $\Phi: \R^d \map \R^n$, that $\Sigma$ is compact, and that the support of $\Phi$ is compact. 
\item
For $(a,P) \in \R^n \times \gr_i$, define the exceptional set
\begin{equation}\label{eap}
E(a,P) := \Big\{ x \in \R^d: \na \Phi(x) \in P, \, \big\langle \Phi(x) -a, \na \Phi(x) \big\rangle = 0 \Big\}.
\end{equation}
Here and throughout, $\langle \bullet,\bullet\rangle$ denotes the Euclidean inner product in $\R^n$, and $\na$ denotes the Euclidean gradient.
\item
Given $(a,P) \in \R^n \times \gr_i$, write $N(a,P)$ for the unique $(n-i)$-dimensional affine plane orthogonal to $P$ and passing through $a$.
\item
In addition, set
\begin{equation}
S_i := \Big\{ (a,P) \in \R^n \times \gr_i : \haus^d\big( E(a,P)\big) > 0 \Big\}.
\end{equation}
\end{itemize}

Now let us prove Theorem \ref{thm: main}. To begin with, we show that the exceptional set $E(a,P)$ is $\haus^d$-null for low dimensional Grassmannians:

\noindent
{\bf Claim 1.} $S_i = \emptyset$ for each $i \leq d$.
\begin{proof}
First we prove for $i<d$. By assumption, the parametrisation map $\Phi: \R^d \map \R^n$ is a $C^1$-embedding, so $\na \Phi$ is of rank $d$ at any point $x \in \R^d$. But for each $x\in E(a,P)$ with $(a,P) \in S_i$, $\na\Phi(x) \in P$, which is an $i$-plane. Thus we get the contradiction.

Now let us consider $i=d$. For $(a,P) \in S_i$ and $x \in E(a,P)$, by construction we have $\Phi(x) \in P \cap N(a,P)$. But $P$ and $N(a,P)$ are clearly transverse at $\Phi(x)$, so $$\dim (P \cap N(a,P)) \leq d-1.$$ That is, $\Phi(x)$ is contained in an affine plane of dimension at most $(d-1)$. However, $\Phi$ is an embedding at $x$, so the inverse function theorem implies that $E(a,P)$ is contained in a submanifold of dimension at most $(d-1)$. This contradicts the definition of $S_i$.  

Hence the assertion follows.  \end{proof}

Next, notice that $\A(\Sigma) = pr_n (S_n)$, as $\p\ball(a,|a-p|)$ is transverse to $\Sigma={\rm image}\,(\Phi)$ at $p=\Phi(x)$ if and only if the line $\Phi(x)-a$ is perpendicular to $T_{\Phi(x)}\Sigma$, which is spanned by $\na\Phi(x)$. One may stratify
\begin{equation}
\A(\Sigma) = pr_n(S_n) \supset pr_{n-1}(S_{n-1}) \supset\cdots\supset pr_d(S_d) = \emptyset,
\end{equation}
where the last equality follows from {\em Claim 1}. Let us set 
\begin{equation}
\A_i := pr_i(S_i) \sim pr_{i-1}(S_{i-1}).
\end{equation}

It remains to show that each $\A_i$ is contained in countably many $(n-d-1)$-affine planes.  For this purpose, we begin with establishing 

\noindent
{\bf Claim 2.} For each $i \in \{d+1, d+2, \ldots, n\}$, the set
\begin{equation}\label{wpi}
\wp_i := \Big\{P \in \gr_i : \text{there exists $a \in \A_i$ such that } (a,P) \in S_i \Big\}
\end{equation}
is at most countable.

\begin{proof}
	First, we show that the sets $E(\bullet,P)$ are essentially disjoint on $\A_i$ for different $P \in \gr_i$. To this end, let us fix $i \in \{d+1, d+2, \ldots, n\}$ and take arbitrary $(a,P)$, $(a', P') \in S_i$ with $a,a'\in\A_i$ and $P\neq P'$ in $\gr_i$. Then, if  $x \in E(a,P) \cap E(a',P')$, one may argue as  for {\em Claim 1} to deduce that	
	\begin{equation*}
	\Phi(x) \in N(a,P) \cap P \cap N(a',P) \cap P'.
	\end{equation*}

The right-hand side is an affine space of dimension at most $(i-1)$, for $P \neq P'$ as $i$-planes. But $a  \in \A_i := pr_i (S_i) \sim pr_{i-1}(S_{i-1})$, so by the definition of $S_{i-1}$ we have  $\haus^{d}\big(E(a,P) \cap E(a',P')\big) = 0$.

Therefore, thanks to the inclusion-exclusion principle, if for $P \neq P'$ in $\gr_i$ and for some $a,a' \in \A_i$ we have both $(a,P) \in S_i$ and $(a',P') \in S_i$, then 
\begin{equation*}
\haus^d\big( E(a,P) \cup E(a',P')\big) = \haus^d \big(E(a,P)\big) + \haus^d \big(E(a',P')\big),
\end{equation*}
and each summand on the right-hand side is positive. On the other hand, all the sets $E(b,Q)$ for $(b,Q) \in S_i$ lie in the support of $\Phi$, which is assumed to be compact.

Hence the assertion follows.  \end{proof}

With {\em Claim 2} at hand, it now remains to show

\noindent
{\bf Claim 3.} Let $i \in \{d+1,d+2,\ldots,n\}$ and $P \in \wp_i$ be fixed. Then the set $\{a\in\A_i : (a,P) \in S_i\}$ is contained in countably many $(n-i)$-affine planes. 

\begin{proof}
	Let us take arbitrary $P \in \wp_i$ and $a \neq \hat{a} \in \A_i$. For any $x\in E(a,P) \cap E(\hat{a},P)$ it then follows that $\Phi(x) \in N(a,P) \cap N(\hat{a}, P)$. However, by elementary Euclidean geometry, we know
	\begin{equation*}
	N(a,P) \cap N(\hat{a}, P) = \begin{cases}
	N(a,P) \equiv N(\hat{a},P)\qquad \text{ if the line through $a,\hat{a}$ is perpendicular to $P$};\\
	\emptyset \qquad \text{ if otherwise}.  
\end{cases}
	\end{equation*}
In particular, 
\begin{equation}\label{Q}
\Big\{a \in \A_i: (a,P) \in S_i\Big\} \subset \bigcup_{Q \in \Xi(i,P)} Q,
\end{equation}
where each $Q$ is an $(n-i)$-dimensional affine plane perpendicular to $P$. 

It is left to prove that $\Xi(i,P)$ can be chosen as a countable set. Indeed, note that for $a \in \A_i$ with $(a,P) \in S_i$, it is covered by $Q\equiv N(a,P)$. Then, \eqref{Q} implies that for $Q \neq \hat{Q}$ in $\Xi(i,P)$ where $Q=N(a,P)$ and $\hat{Q}=N(\hat{a},P)$, one has 
\begin{equation*}
\haus^d \big(E(a,P) \cap E(\hat{a}, P)\big) = 0.
\end{equation*} 
The same inclusion-exclusion principle argument as in the proof of {\em Claim 2} yields that there are at most countably many such $Q$'s.

Hence the assertion follows. \end{proof}

To summarise, we have proved that
\begin{equation}
\A(\Sigma) = \bigsqcup_{i=d+1}^n \bigcup_{P \in \wp_i}\bigcup_{Q \in \Xi(i,P)} Q.
\end{equation}
Here each $Q$ is an $(n-i)$-dimensional affine space with $d+1 \leq i \leq n$; also, $\wp_i$ and $\Xi(i,P)$ are countable indexing sets. Therefore, the proof of Theorem \ref{thm: main} is now complete.

\begin{remark}
The assumption $d \geq 1$ in Theorem \ref{thm: main} is necessary. Indeed, if $\Sigma$ is a non-empty $0$-dimensional manifold, {\it i.e.}, a discrete set with the counting measure, then it is not transverse to any $(n-1)$-dimensional manifold. It follows that $\mathscr{A}(\Sigma)=\R^n$, which cannot be contained in a countable union of $(n-1)$-dimensional affine planes.
\end{remark}

\section{An Example}\label{sharp}

In this section, we construct an explicit example of a $d$-dimensional $C^1$-submanifold $\Sigma \subset \R^n$, such that  $\A(\Sigma)$ is exactly equal to a countable collection of $(n-d-1)$-dimensional affine planes. Thus, the rectifiability dimension $(n-d-1)$ is sharp.

\medskip
{\bf 1.} We first describe the example for the special case $n=3$ and $d=1$, as this is easy to be pictured in $\R^3$.

 Consider
\begin{equation*}
\Sigma_0 := \bigsqcup_{i\in\mathbb{Z}}\bigg\{ \Big( \frac{1}{4}\cos\theta, \frac{1}{4}\sin\theta, 0 \Big) + (i,0,0) : \theta \in [0,2\pi] \bigg\}.
\end{equation*}
It is a disjoint union of countably many $\mathbb{S}^1$ in the $(x,y)$-plane with $z$-coordinate equal to $0$. Then we have 
\begin{equation*}
\A(\Sigma_0) = \bigsqcup_{i\in\mathbb{Z}} \Big\{ \text{the $z$-axis} + (i,0,0)\Big\}.
\end{equation*}
Indeed, let us check that $\A(\gamma)=\text{the $z$-axis}$, where $\gamma:=\{(4^{-1}\cos\theta, 4^{-1}\sin\theta,0): \theta \in [0,2\pi]\}$; then it follows by translation. For each point $a$ on the $z$-axis and each $p \in \gamma$, $\gamma$ is a meridian of the $2$-sphere $\p\ball(a, |a-p|)$. Thus $T_p\gamma \subset T_p \p\ball(a, |a-p|)$, so $\gamma$ and $\p\ball(a, |a-p|)$ are not transverse at $p$. This shows the $z$-axis $\subset \A(\gamma)$. In fact, it is clear that the circle $\gamma$ is not transverse to a $2$-sphere $S$ in $\R^3$ at a set of positive $\haus^1$-measure if and only if $\gamma$ is a meridian of $S$. The other inclusion thus follows.

Moreover, we can modify the example $\Sigma_0$ to make it a connected manifold. For example, we can ``open up'' $\Sigma_0$ at the ends of each circle, connect the adjacent components, and then make it smooth by mollification. More precisely, we define $\Sigma_2$ as follows. 

First, take $\e \leq 100^{-1}$ and set \begin{equation*}
\Sigma_1 := \bigsqcup_{i\in\mathbb{Z}}\bigg\{ \Big( \frac{1}{4}\cos\theta, \frac{1}{4}\sin\theta, 0 \Big) + (i,0,0) : \theta \in [2\e, \pi -2\e] \sqcup [\pi+2\e, 2\pi - 2\e] \bigg\}.
\end{equation*}
Then, for each connected component of $\Sigma_1$, we join the following pairs of points by a straight line segment for each $i \in \mathbb{Z}$ (these points thus become ``corners''): 
\begin{equation*}
\begin{cases}
\Big( \frac{1}{4}\cos 2\e, \frac{1}{4}\sin 2\e, 0 \Big) + (i,0,0)\quad\text{and}\quad\Big( \frac{1}{4}\cos(\pi -2\e), \frac{1}{4}\sin (\pi -2\e), 0 \Big) + (i+1,0,0);  \\
\Big( \frac{1}{4}\cos (2\pi -2 \e), \frac{1}{4}\sin (2\pi - 2\e), 0 \Big) + (i,0,0)\quad\text{and}\quad\Big( \frac{1}{4}\cos(\pi +2\e), \frac{1}{4}\sin (\pi + 2\e), 0 \Big) + (i+1,0,0).
\end{cases}
\end{equation*}
Finally, we smooth up the corners to get $\Sigma_2$: one may ensure that the set is unchanged outside $\e$-neighbourhoods of the corners. By definition, we still have $\A(\Sigma_2) = $ the $z$-axis and its translates by $(i,0,0)$, $i \in \mathbb{Z}$. Here $\Sigma_2$ is a connected $1$-dimensional $C^\infty$-submanifold of $\R^3$, whose $\A(\Sigma_2)$ equals to a countable collection of $1$-dimensional affine planes.

In effect, we have constructed 
$$
\Sigma_2 = \mathbb{S}^1 \, \# \, \mathbb{S}^1 \, \# \, \mathbb{S}^1 \, \# \, \cdots
$$
topologically, where each connected summand $\mathbb{S}^1$ sits in $\R^2_i  \times \{0\} \subset \R^3$, $i \in \mathbb{Z}$ and $\R^2_i =\R^2 + (i,0,0)$. $\mathbb{S}^1$ denotes the topological $1$-sphere (suitably scaled), and the connected sum operation $\#$ is taken by using very thin necks ($\simeq 3\e$, upon mollification).  Then, $\A(\Sigma_2)$ consists of the translated copies of $\R$, which is the orthogonal complement of  $\R^2_i  \times \{0\}$ in $\R^3$.

\medskip
{\bf 2.} We can generalise the above constructions to arbitrary $n$ and $d$.

Extending the idea at the end of {\bf 1.}, let us construct $\Sigma$ by $$
\Sigma = \mathbb{S}^d \, \# \, \mathbb{S}^d \, \# \, \mathbb{S}^d \, \# \, \cdots.
$$
Here, each $\mathbb{S}^d$ is a round sphere ({\it e.g.}, scaled by $4^{-n}$) embedded in $\R^{d+1}$, and we view $\R^{d+1}$ as consisting of the first $(d+1)$-coordinates in $\R^n$. The connected sum operation again uses arbitrarily thin necks. The spheres $\mathbb{S}^d$ are congruent to each other and obtained via translating the one centred at the origin (call it $\mathbb{S}^d_0$) by $(i,0,\ldots,0)$, where $i \in \mathbb{Z}$. 

Take an arbitrary point
\begin{equation*}
a \in \R^{n-d-1} + (i,0,\ldots,0);
\end{equation*}
$\R^{n-d-1}$ is the orthogonal complement vector subspace of $\R^n$ to $\R^{d+1}$, {\it i.e.} containing the last $(n-d-1)$ coordinates in $\R^n$. Also, take an arbitrary point 
\begin{equation*}
p \in \mathbb{S}^d_0 + (i,0,\ldots,0) \equiv \mathbb{S}^d_i.
\end{equation*}
Then $\p\ball(a,|a-p|)$ is the $(n-1)$-sphere centred at $a$ with radius $|a-p|$. But $\mathbb{S}^d_i \subset \p\ball(a,|a-p|)$ and $T_p\mathbb{S}^d_i \subset T_p\p\ball(a,|a-p|)$, as $\mathbb{S}^d_i$ is an (iterated) meridian of $\p\ball(a,|a-p|)$. Together with the definition of $\A(\Sigma)$ and the thinness of the necks, this implies that 
\begin{equation*}
 \R^{n-d-1} + (i,0,\ldots,0) \subset \A(\Sigma)\qquad \text{ for each } i \in\mathbb{Z}.
\end{equation*}

On the other hand, if $a \notin  \R^{n-d-1} + (i,0,\ldots,0)$ for some $i\in\mathbb{Z}$, then for all but finitely many $p \in \mathbb{S}^i_d$, there is a neighbourhood $U$ of $p$ on $\mathbb{S}^i_d$ such that $|a-q|$ for $q\in U$ is not extremised at $p$. Thus $\p\ball(a,|a-p|)$ is transverse to $\mathbb{S}^d_i$ (hence to $\Sigma$) at $p$. This implies $a\notin \A(\Sigma)$, {\it i.e.}, 
\begin{equation*}
 \R^{n-d-1} + (i,0,\ldots,0) \supset \A(\Sigma)\qquad \text{ for each } i \in\mathbb{Z}.
\end{equation*}

Thus, we have constructed $\Sigma$ so that $\A(\Sigma)$ equals to a countable collection of $(n-d-1)$-dimensional affine planes.

\bigskip
\noindent
{\bf Acknowledgement}.
This work has been done during Siran Li's stay as a CRM--ISM postdoctoral fellow at  the Centre de Recherches Math\'{e}matiques, Universit\'{e} de Montr\'{e}al and the Institut des Sciences Math\'{e}matiques. Siran Li would like to thank these institutions for their hospitality. The author also thanks the referee for pointing out a mistake in an earlier version of the draft.

\end{document}